\documentclass[12pt]{article}
\usepackage{mathrsfs}
\usepackage{epic,eepic,epsf,epsfig}
\usepackage{amsfonts,srcltx,mathrsfs}
\usepackage{multirow}
\usepackage{amsmath}
\usepackage{amssymb}
\usepackage{amsbsy}
\usepackage{graphicx}
\usepackage{amsfonts}
\usepackage{color}
\usepackage{setspace}

\newtheorem{lem}{Lemma}[section]%
\newtheorem{theorem}[lem]{Theorem}%

\def\nd{\mathrel{\bigm|\kern-.7em/}}

\def\f{\noindent}

\def\P\GammaL{\hbox{\rm P\GammaL}}

\begin{document}
\title{A spectral condition for Hamilton cycles in tough bipartite graphs}

\footnotetext{* Corresponding author}
\footnotetext{E-mails: sdlgaly@126.com, zhangwq@pku.edu.cn}

\author{Lianyang Ai, Wenqian Zhang*\\
{\small School of Mathematics and Statistics, Shandong University of Technology}\\
{\small Zibo, Shandong 255000, P.R. China}}
\date{}
\maketitle

\begin{abstract}
Let $G$ be a graph. The {\em spectral radius} of $G$ is the largest eigenvalue of its adjacency matrix. For a non-complete bipartite graph $G$ with parts $X$ and $Y$, the {\em bipartite toughness} of $G$ is defined as $t^{B}(G)=\min\left\{\frac{|S|}{c(G-S)}\right\}$,
where the minimum is taken over all proper subsets $S\subset X$ (or $S\subset Y$) such that $c(G-S)>1$.
In this paper, we give a sharp spectral radius condition for balanced bipartite graphs $G$ with $t^{B}(G)\geq1$ to guarantee that $G$ contains Hamilton cycles. This solves a problem proposed in \cite{CFL}.
\bigskip

\f {\bf Keywords:} spectral radius; Hamilton cycle; bipartite graph; $[a,b]$-factor.

\end{abstract}

\baselineskip 17 pt

\section{Introduction}

All graphs considered here are finite, undirected and simple. For a graph $G$, let $V(G)$ and $E(G)$ denote the vertex set and edge set of $G$, respectively. Let $e(G)=|E(G)|$. For a vertex $u$, let $d_{G}(u)$ denote its degree. Let $\delta(G)$ be the minimum degree of $G$. For a subset $S$ of $V(G)$, let $G[S]$ be the subgraph induced by $S$, and let $G-S$ obtained from $G$ by deleting the vertices in $S$ and the incident edges. For two disjoint subsets $U$ and $V$ of $V(G)$, let $e_{G}(U,V)$ be the number of edges between $U$ and $V$. For certain integers $m,n$, let $K_{n}$ and $K_{m,n}$ denote the complete graph and the complete bipartite graph, respectively. The {\em spectral radius} of $G$ is the largest eigenvalue of its adjacency matrix. For any terminology used but not defined here, one may refer to \cite{BH}.

Let $G$ be a graph. For non-negative integers $a$ and $b$, an $[a,b]$-factor of $G$ is a spanning subgraph with degrees between $a$ and $b$. An $[a,b]$-factor is called a $k$-factor if $a=b=k$. Recently, many researchers studied the spectral radius conditions for graphs to contain $[a,b]$-factors (for example, see \cite{CS,CFL,FL,FL2,FLL2,HL,HLY,O,WS}). A Hamilton cycle of $G$ is a connected $2$-factor. Recently, many researchers studied the spectral radius conditions for graphs to contain Hamilton cycles (for example, see \cite{Bene,FV,FLL,LLT,LSX,LN,NG,ZBWL,ZBWL2}).

Chv\'{a}tal \cite{C} defined the toughness $t(G)$ of a non-complete graph $G$ as $t(G)=\min\left\{\frac{|S|}{c(G-S)}\right\}$,
where the minimum is taken over all proper subsets $S\subset V(G)$, and $c(G-S)$ denotes the number of components of $G-S$. We say that $G$ is $k$-tough, if $t(G)\geq k$. Clearly, if $G$ has a Hamiltonian cycle then $G$ is 1-tough. Conversely, Chv\'{a}tal \cite{C}
proposed the following conjecture.

\medskip

\f{\bf Conjecture.} {\rm (\cite{C})} There exists a constant $C$ such that every $C$-tough graph is Hamiltonian.

\medskip

Bauer, Broersma and Veldman \cite{BBV} constructed examples which implies that $C\geq\frac{9}{4}$ (if it exists). In general, this conjecture still remains open.

For $n\geq8$, denote the vertex set of $K_{n-3}$ by $\left\{v_{1},v_{2},...,v_{n-3}\right\}$. Let $M_{n}$ be the graph obtained from $K_{n-3}$ by adding 3 vertices $u_{1},u_{2},u_{3}$ (all) connecting to $v_{1}$, and connecting $u_{i}$ to $v_{i+1}$ for any $1\leq i\leq3$. Clearly, $M_{n}$ is a 1-tough graph of order $n$. Note that $M_{n}$ contains no Hamilton cycles, since the three vertices $u_{1},u_{2},u_{3}$ are of degree 2 and have a common neighbor $v_{1}$. Fan, Lin and Lu \cite{FLL} proposed the problem: What is the spectral condition to guarantee the existence of a Hamiltonian cycle among $t$-tough graphs? They proved the following result.

\begin{theorem}
Suppose that $G$ is a connected 1-tough graph of order $n\geq18$ (with $\delta(G)\geq2$). If
$\rho(G)\geq\rho(M_{n})$, then $G$ contains a Hamiltonian cycle, unless $G=M_{n}$.
\end{theorem}

Let $G=G(X,Y)$ denote a bipartite graph with parts $X$ and $Y$, where $|X|\leq |Y|$. The graph $G$ is balanced if $|X|=|Y|$.  Observe that $t(G)\leq1$, since $G-X$ has $|Y|$ components. Bian \cite{B} introduced the concept of bipartite toughness for $G$. For a non-complete bipartite graph $G$, its bipartite toughness $t^{B}(G)$ is defined
as $t^{B}(G)=\min\left\{\frac{|S|}{c(G-S)}\right\}$,
where the minimum is taken over all proper subsets $S\subset X$ (or $S\subset Y$) such that $c(G-S)>1$.  We say that $G$ is $k$-tough, if $t^{B}(G)\geq k$.

For $n\geq5$, denote the two parts of $K_{n,n-3}$ by $\left\{u_{1},u_{2},...,u_{n}\right\}$ and $\left\{v_{4},v_{5},...,v_{n}\right\}$. Let $G_{n,n}$ be the graph obtained from $K_{n,n-3}$ by adding 3 vertices $v_{1},v_{2},v_{3}$ (all) connecting to $u_{1}$, and connecting $v_{i}$ to $u_{i+1}$ for any $1\leq i\leq3$. Clearly, $G_{n,n}$ is a 1-tough balanced bipartite graph of order $2n$. Note that $G_{n,n}$ contains no 2-factors, since the three vertices $v_{1},v_{2},v_{3}$ are of degree 2 and have a common neighbor $u_{1}$.
Chen, Fan and Lin \cite{CFL} proved the following result.

\begin{theorem}{\rm (\cite{CFL})} \label{2-factor}
Suppose that $G$ is a connected 1-tough balanced bipartite graph of order $2n$. If
$\rho(G)\geq\rho(G_{n,n})$, then $G$ contains a 2-factor, unless $G=G_{n,n}$.
\end{theorem}

At the end of \cite{CFL}, the authors guessed that $G$ contains a Hamilton cycle in the condition of Theorem \ref{2-factor}. In this paper, we confirm this.

\begin{theorem}\label{main1}
Let $G=G(X,Y)$ be a 1-tough balanced bipartite graph of order $2n$, where $|X|=|Y|=n\geq16$. If $\rho(G)\geq\rho(G_{n,n})$, then $G$ contains a Hamilton cycle, unless $G=G_{n,n}$.
\end{theorem}

\section{Proof of Theorem \ref{main1} }

To prove Theorem \ref{main1}, we first include some lemmas.

\begin{lem}{\rm(\cite{BH})}\label{subgraph}
Let $G$ be a connected graph. For any subgraph $H$ of $G$, $\rho(H)\leq\rho(G)$ with equality if and only if $H=G$.
\end{lem}

\begin{lem}{\rm(\cite{BFP})}\label{bi-graph}
If $G$ is a bipartite graph, then $\rho(G)\leq\sqrt{e(G)}$. Moreover, equality holds if and only if $G$ is a complete bipartite graph with some isolated vertices.
\end{lem}

Let $G=G(X,Y)$ be a balanced bipartite graph of order $2n$. The {\em bipartite closure} $H$  of $G$, is the graph obtained from $G$ by recursively joining
pairs of non-adjacent vertices $u\in X$ and $v\in Y$ satisfying $d_{G}(u)+d_{G}(v)\geq n+1$ until no such pair remains. The following is the bipartite closure technique, which is widely used in the study of Hamilton cycles in balanced bipartite graphs.

\begin{theorem}{\rm(\cite{BC})}\label{bi-closure}
 A balanced bipartite graph $G$ is Hamiltonian if and only if its  bipartite closure is Hamiltonian.
\end{theorem}

Now we begin the proof of Theorem \ref{main1}.

\medskip

\f{\bf Proof of Theorem \ref{main1}:} Since $G$ is 1-tough, we see $\delta(G)\geq2$. Let $H$ be the bipartite closure of $G$. Then $H$ is a 1-tough balanced bipartite graph.  Note that $G$ is a spanning subgraph of $H$. By Lemma \ref{subgraph}, $\rho(H)\geq\rho(G)$, and equality holds if and only if $G=H$. By Theorem \ref{bi-closure}, if $H$ contains a Hamilton cycle, then $G$ also contains one. Thus, we can assume that $H$ contains no Hamilton cycles.  It suffices to show $H=G_{n,n}$ (implying $G=G_{n,n}$ as $\rho(G)\geq\rho(G_{n,n})$). Since $G_{n,n}$ contains $K_{n,n-3}$ as a subgraph, by Lemma \ref{subgraph} we have $\rho(H)>\rho(K_{n,n-3})=\sqrt{n(n-3)}$. By Lemma \ref{bi-graph}, $\rho(H)\leq\sqrt{e(H)}$. Hence,
$$e(H)>n(n-3)=n^{2}-3n.$$

Recall that $H$ is the bipartite closure of $G$. Thus, $u\in X$ is adjacent to $v\in Y$ whenever $d_{H}(u)+d_{H}(v)\geq n+1$. Let $X_{0}\subseteq X$ and $Y_{0}\subseteq Y$ be the sets of vertices of degree at least $\frac{n+1}{2}$ in $H$. Then all the vertices in $X_{0}$ are adjacent to all the vertices in $Y_{0}$.

Now we show that $|X_{0}|\geq\frac{n+1}{2}$ and $|Y_{0}|\geq\frac{n+1}{2}$. Suppose not (for example, $|X_{0}|\leq\frac{n}{2}$). Note that $d_{H}(u)\leq\frac{n}{2}$ for any $u\in X-X_{0}$. Thus, noting $n\geq16$,
\begin{equation}
\begin{aligned}
e(H)&=\sum_{v\in X}d_{H}(v)\\
&\leq n|X_{0}|+\frac{n}{2}(n-|X_{0}|)\\
&\leq\frac{3}{4}n^{2}\\
&\leq n^{2}-3n,
\end{aligned}\notag
\end{equation}
contradicting the fact that $e(H)>n^{2}-3n$. Hence, $|X_{0}|\geq\frac{n+1}{2}$ and $|Y_{0}|\geq\frac{n+1}{2}$.

Set $s=|X_{0}|$ and $ t=|Y_{0}|$. Without loss of generality, assume that $s\geq t$.

\medskip

\f{\bf Claim 1.} For any $u\in X-X_{0}$ and $v\in Y-Y_{0}$, $d_{H}(u)\leq n-s$ and $d_{H}(v)\leq n-t$.

\medskip

\f{\bf Proof of Claim 1.} Let $u\in X-X_{0}$. Noting that $d_{H}(u)\leq\frac{n}{2}$ and $t\geq\frac{n+1}{2}$, there is some $v'\in Y_{0}$ which is not adjacent to $u$. Note that $d_{H}(v')\geq s$ as $v'$ is adjacent to all the vertices in $X_{0}$. Thus, $d_{H}(u)\leq n-d_{H}(v')\leq n-s$. Similarly, we can show that $d_{H}(v)\leq n-t$ for any $v\in Y-Y_{0}$. \hfill$\Box$

\medskip

\f{\bf Claim 2.} $s,t\geq n-3$ and $s,t\neq n-1$.

\medskip

\f{\bf Proof of Claim 2.} We first show $s,t\geq n-3$. Suppose not. Then $t\leq n-4$. By Claim 1, $d_{H}(v)\leq n-t$ for any $v\in Y-Y_{0}$. Thus,

 \begin{equation}
\begin{aligned}
e(H)&=\sum_{v\in Y}d_{H}(v)\\
&=\sum_{v\in Y_{0}}d_{H}(v)+\sum_{v\in Y-Y_{0}}d_{H}(v)\\
&\leq nt+(n-t)(n-t)\\
&=n^{2}-t(n-t)\\
&\leq n^{2}-4(n-4)\\
&\leq n^{2}-3n,
\end{aligned}\notag
\end{equation}
a contradiction. Hence, $s,t\geq n-3$.

Now we show that $s,t\neq n-1$. Suppose not (for example, $s= n-1$). Let $u$ be the only vertex in $X-X_{0}$. By Claim 1, $d_{H}(u)\leq n-(n-1)=1$, contradicting the fact that $\delta(H)\geq2$. Thus $s\neq n-1$. Similarly, we have $t\neq n-1$. \hfill$\Box$

 If $t=n$, then $s=n$. Thus, $H$ is the complete graph $K_{n,n}$, which contains a Hamilton cycle. This contradicts the assumption that $H$ contains no Hamilton cycles. Thus $t=n-2$ or $t=n-3$ by Claim 2. Denote by $X=\left\{u_{1},u_{2},...,u_{n}\right\}$ and $Y=\left\{v_{1},v_{2},...,v_{n}\right\}.$ Without loss of generality, assume that $X-X_{0}=\left\{u_{1},u_{2},\ldots,u_{n-s}\right\}$ and $Y-Y_{0}=\left\{v_{1},v_{2},\ldots,v_{n-t}\right\}.$ Now we prove the theorem by two cases upon the values of $t$.

\medskip

\f{\bf Case 1.} $t=n-2$.

\medskip

Then $s=n-2$ or $s=n$ by Claim 2.

\medskip

\f{\bf Subcase 1.1.} $s=n-2$.

\medskip

In this case, we have $X-X_{0}=\left\{u_{1},u_{2}\right\}$ and $Y-Y_{0}=\left\{v_{1},v_{2}\right\}.$ For any $1\leq i\leq 2$, we have $d_{H}(u_{i})\leq2$ and $d_{H}(v_{i})\leq2$ by Claim 1. Hence, $d_{H}(u_{i})=2$ and $d_{H}(v_{i})=2$ as $\delta(H)\geq2$.
Then,
 \begin{equation}
\begin{aligned}
e(H)&\leq e_{H}(X_{0},Y_{0})+\sum_{1\leq i\leq 2}(d_{H}(u_{i})+d_{H}(v_{i}))\\
&\leq (n-2)^{2}+2\cdot 4\\
&\leq n^{2}-3n,
\end{aligned}\notag
\end{equation}
a contradiction.

\medskip

\f{\bf Subcase 1.2.} $s=n$.

\medskip

In this case, $X_{0}=X$ and $Y-Y_{0}=\left\{v_{1},v_{2}\right\}$. Recall that $H[X\cup Y_{0}]=K_{n,n-2}$. Similar to Subcase 1.1, we can show that $d_{H}(v_{i})=2$ for any $1\leq i\leq2$. Without loss of generality, assume that $v_{1}$ is adjacent to $u_{1}$ and $u_{2}$. Since $H$ is 1-tough, $v_{2}$ has a neighbor distinct from $u_{1}$ and $u_{2}$, say $u_{3}$. If $v_{1},v_{2}$ have a common neighbor, say $u_{1}$, then $u_{2}v_{1}u_{1}v_{2}u_{3}v_{3}u_{4}v_{4}\cdots u_{n}v_{n}u_{2}$ is a Hamilton cycle, a contradiction. If $v_{1},v_{2}$ have no common neighbors, then $v_{2}$ has another neighbor, say $u_{4}$. Then $u_{1}v_{1}u_{2}v_{3}u_{3}v_{2}u_{4}v_{4}u_{5}v_{5}\cdots u_{n}v_{n}u_{1}$ is a Hamilton cycle, a contradiction.

\medskip

\f{\bf Case 2.} $t=n-3$.

\medskip

\f{\bf Subcase 2.1.} $s=n-3$.

\medskip

By Claim 1, $d_{H}(u_{i})\leq3$ and $d_{H}(v_{i})\leq3$ for any $1\leq i\leq 3$.
 Then,
 \begin{equation}
\begin{aligned}
e(H)&\leq e_{H}(X_{0},Y_{0})+\sum_{1\leq i\leq 3}(d_{H}(u_{i})+d_{H}(v_{i}))\\
&\leq (n-3)^{2}+3\cdot 6\\
&\leq n^{2}-3n,
\end{aligned}\notag
\end{equation}
a contradiction.

\medskip

\f{\bf Subcase 2.2.} $s=n-2$.

\medskip

In this case, $X-X_{0}=\left\{u_{1},u_{2}\right\}$ and $Y-Y_{0}=\left\{v_{1},v_{2},v_{3}\right\}.$ By Claim 1, $d_{H}(u_{i})\leq2$ for $1\leq i\leq2$, and $d_{H}(v_{j})\leq3$ for any $1\leq j\leq 3$. Then
\begin{equation}
\begin{aligned}
e(H)&\leq e_{H}(X_{0},Y_{0})+\sum_{1\leq i\leq 2}d_{H}(u_{i})+\sum_{1\leq j\leq 3}d_{H}(v_{i})\\
&\leq (n-3)(n-2)+2\cdot2+3\cdot 3\\
&\leq n^{2}-3n,
\end{aligned}\notag
\end{equation}
a contradiction.

\medskip

\f{\bf Subcase 2.3.} $s=n$.

\medskip

In this case, $X_{0}=X$ and  $Y-Y_{0}=\left\{v_{1},v_{2},v_{3}\right\}$. Recall that $H[X\cup Y_{0}]=K_{n,n-3}$. By Claim 1, $d_{H}(v_{i})\leq3$ for any $1\leq i\leq3$.

A {\em good linear forest} in $H$ is the disjoint union of pathes $\cup_{ i\geq1 }Q_{i}$, where each path $Q_{i}$ has its two end-vertices in $X$ and $V(Q_{i})\cap Y\subseteq\left\{v_{1},v_{2},v_{3}\right\}$. Clearly, a good linear forest has at most 3 disjoint paths.

\medskip

\f{\bf Claim 3.} $H$ has no good linear forests.

\medskip

\f{\bf Proof of Claim 3.} Suppose not. Let $\cup_{ 1\leq i\leq i_{0} }Q_{i}$ be a good linear forest in $H$. If $i_{0}=1$, without loss of generality, assume that $Q_{1}=u_{1}v_{1}u_{2}v_{2}u_{3}v_{3}u_{4}$.
Then  $$u_{1}v_{1}u_{2}v_{2}u_{3}v_{3}u_{4}v_{4}u_{5}v_{5}\cdots u_{n}v_{n}u_{1}$$
 is a Hamilton cycle, a contradiction.  If $i_{0}=2$, without loss of generality, assume that $Q_{1}=u_{1}v_{1}u_{2}v_{2}u_{3}$ and $Q_{2}=u_{4}v_{3}u_{5}$.
Then  $$u_{1}v_{1}u_{2}v_{2}u_{3}v_{4}u_{4}v_{3}u_{5}v_{5}u_{6}v_{6}\cdots u_{n}v_{n}u_{1}$$
 is a Hamilton cycle, a contradiction.
 If $i_{0}=3$, without loss of generality, assume that $Q_{1}=u_{1}v_{1}u_{2},Q_{2}=u_{3}v_{2}u_{4}$ and $Q_{2}=u_{5}v_{3}u_{6}$.
Then  $$u_{1}v_{1}u_{2}v_{4}u_{3}v_{2}u_{4}v_{5}u_{5}v_{3}u_{6}v_{6}u_{7}v_{7}\cdots u_{n}v_{n}u_{1}$$
 is a Hamilton cycle, a contradiction. This finishes the proof of Claim 3. \hfill$\Box$

\medskip

\f{\bf Claim 4.} $d_{H}(v_{i})=2$ for any $1\leq i\leq3$.

\medskip

\f{\bf Proof of Claim 4.} Suppose not. Without loss of generality, assume $d_{H}(v_{1})=3$, and $v_{1}$ is adjacent to $u_{1},u_{2},u_{3}$. We first consider the case that one of $\left\{v_{2},v_{3}\right\}$, say $v_{2}$, is only adjacent to the vertices in  $\left\{u_{1},u_{2},u_{3}\right\}$. Without loss of generality, we can assume that $v_{2}$ is adjacent to $u_{1},u_{2}$ (and possibly $u_{3}$). Since $H$ is 1-tough, $v_{3}$ has a neighbor in $X-\left\{u_{1},u_{2},u_{3}\right\}$, say $u_{4}$, and has another neighbor $u_{\ell}$ with $\ell\neq4$. If $1\leq\ell\leq3$, we only consider the case of $\ell=1$ (as it is similar for other cases). That is to say, $v_{3}$ is adjacent to $u_{1}$. Then  $u_{3}v_{1}u_{2}v_{2}u_{1}v_{3}u_{4}$ is a good linear forest, which is a contradiction by Claim 3.
If $\ell\geq5$, without loss of generality, assume that $\ell=5$. Then the union of $u_{1}v_{2}u_{2}v_{1}u_{3}$ and $u_{4}v_{3}u_{5}$ is a good linear forest, a contradiction.

It remains the case that both $v_{2}$ and $v_{3}$ are adjacent to some vertices in  $X-\left\{u_{1},u_{2},u_{3}\right\}$. We only assume that both $v_{2}$ and $v_{3}$ have only one neighbor in  $X-\left\{u_{1},u_{2},u_{3}\right\}$. (It is more easy to find a good linear forest if $v_{2}$ or $v_{3}$ has distinct neighbors in  $X-\left\{u_{1},u_{2},u_{3}\right\}$.) First consider that the only neighbors of $v_{2}$ and $v_{3}$ in  $X-\left\{u_{1},u_{2},u_{3}\right\}$ are the same, say $u_{4}$. Since $H$ is 1-tough, $v_{2}$ and $v_{3}$ must have distinct neighbors in  $\left\{u_{1},u_{2},u_{3}\right\}$. Without loss of generality, we can assume that $v_{2}$ is adjacent to $u_{1}$, and $v_{3}$ is adjacent to $u_{2}$.  Then  $u_{1}v_{2}u_{4}v_{3}u_{2}v_{1}u_{3}$ is a good linear forest, a contradiction. It remains that the only neighbors of $v_{2}$ and $v_{3}$ in  $X-\left\{u_{1},u_{2},u_{3}\right\}$ are distinct. Without loss of generality, we can assume that $v_{2}$ is adjacent to $u_{4}$, and $v_{3}$ is adjacent to $u_{5}$. Since $\delta(H)\geq2$, we have that  $v_{2}$ is adjacent to $u_{i}$, and $v_{3}$ is adjacent to $u_{j}$, where $1\leq i,j\leq 3$. If $u_{i}=u_{j}$, without loss of generality, assume that both $v_{2}$ and $v_{3}$ are adjacent to $u_{1}$. Then  the union of $u_{2}v_{1}u_{3}$ and $u_{4}v_{2}u_{1}v_{3}u_{5}$ is a good linear forest, a contradiction.
If $u_{i}\neq u_{j}$, without loss of generality, assume that $v_{2}$ is adjacent to $u_{1}$, and $v_{3}$ is adjacent to $u_{2}$. Then  the union of $u_{1}v_{2}u_{4}$ and $u_{3}v_{1}u_{2}v_{3}u_{5}$ is a good linear forest, a contradiction.
This finishes the proof of Claim 4. \hfill$\Box$

\medskip

By Claim 4, $d_{H}(v_{i})=2$ for any $1\leq i\leq3$. Now consider the subgraph $H_{0}$ of $H$ induced by the edges incident with some $v_{i}$ for $1\leq i\leq3$. Since $H$ is 1-tough, we see $7\leq|H_{0}|\leq9$. If the maximum degree of $H_{0}$ is 2, then (noting that $H$ is 1-tough), $H_{0}$ must be a good linear forest, a contradiction by Claim 3. Therefore, the maximum degree of $H_{0}$ must be 3. That is to say, there is some vertex in $X$, say $u_{1}$, which is adjacent to all the three vertices $v_{1},v_{2}$ and $v_{3}$. Since $H$ is 1-tough, the other neighbors (except $u_{1}$) of $v_{1},v_{2}$ and $v_{3}$ must be distinct. Hence, $H=G(n,n)$.
This completes the proof. \hfill$\Box$

\medskip

\f{\bf Data availability statement}

\medskip

There is no associated data.

\medskip

\f{\bf Declaration of Interest Statement}

\medskip

There is no conflict of interest.

\end{document}